\newcommand{\R}{\mathds R}
\newcommand{\Z}{\mathds Z}
\newcommand{\Spl}{\mathrm{Sp}}
\newcommand{\GL}{\mathrm{GL}}
\newcommand{\Dim}{\mathrm{dim}}
\newcommand{\codim}{\mathrm{codim}}
\newcommand{\Bsym}{\mathrm{B_{\mathrm{sym}}}}
\newcommand{\segnatura}{\mathrm{sign}}

\newcommand{\Ker}{\mathrm{Ker}}

\newcommand{\Ddt}{\tfrac{\mathrm D}{\mathrm dt}}

\documentclass[twoside,final,a4paper,10pt]{amsart}

\usepackage{times}
\usepackage{dsfont}
\usepackage[all]{xy}

\renewcommand{\contentsline}[3]{\csname new#1\endcsname{#2}{#3}}
\newcommand{\newchapter}[2]{\bigskip\hbox to \hsize{\vbox{\advance\hsize by -.5cm\baselineskip=12pt\parfillskip=0pt\leftskip=2cm\noindent\hskip -2cm #1\leaders\hbox{.}\hfil\hfil\par}$\,$#2\hfil}}
\newcommand{\newsection}[2]{\medskip\hbox to \hsize{\vbox{\advance\hsize by -.5cm\baselineskip=12pt\parfillskip=0pt\leftskip=2.5cm\noindent\hskip -2cm #1\leaders\hbox{.}\hfil\hfil\par}$\,$#2\hfil}}
\newcommand{\newsubsection}[2]{\medskip\hbox to \hsize{\vbox{\advance\hsize by -.5cm\baselineskip=12pt\parfillskip=0pt\leftskip=3.5cm\noindent\hskip -2cm #1\leaders\hbox{.}\hfil\hfil\par}$\,$#2\hfil}}

\numberwithin{equation}{section}

\title[Comparison results for conjugate and focal points]{Comparison results for conjugate and focal points in semi-Riemannian geometry via Maslov index}
\author[M. \'A. Javaloyes]{Miguel \'Angel Javaloyes}
\address{Departamento de Geometr\'{\i}a y Topolog\'{\i}a.\hfill\break\indent
 Facultad de Ciencias, Universidad de Granada.
 \hfill\break\indent
 Campus Fuentenueva s/n, 18071 Granada, Spain}
 \email{ma.javaloyes@gmail.com}
\author[P.\ Piccione]{Paolo Piccione}
\address{Departamento de Matem\'atica,\hfill\break\indent
Universidade de S\~ao Paulo, \hfill\break\indent Rua do Mat\~ao
1010,\hfill\break\indent CEP 05508-900, S\~ao Paulo, SP, Brazil}
\email{piccione.p@gmail.com}
\thanks{The first author is partially supported by Regional J.
Andaluc\'{\i}a Grant P06-FQM-01951 and by Spanish MEC Grant MTM2007-64504. The second author is partially sponsored by CNPq and Fapesp, Brazil.}


\date{July 31st, 2008}

\begin{document}

% Theorems and such

\theoremstyle{plain}\newtheorem*{teon}{Theorem}
\theoremstyle{definition}\newtheorem*{defin*}{Definition}
\theoremstyle{plain}\newtheorem{teo}{Theorem}[section]
\theoremstyle{plain}\newtheorem{prop}[teo]{Proposition}
\theoremstyle{plain}\newtheorem{lem}[teo]{Lemma}
\theoremstyle{plain}\newtheorem{cor}[teo]{Corollary}
\theoremstyle{definition}\newtheorem{defin}[teo]{Definition}
\theoremstyle{remark}\newtheorem{rem}[teo]{Remark}
\theoremstyle{plain} \newtheorem{assum}[teo]{Assumption}
\swapnumbers
\theoremstyle{definition}\newtheorem{example}{Example}[section]
\theoremstyle{plain} \newtheorem*{acknowledgement}{Acknowledgements}
\theoremstyle{definition}\newtheorem*{notation}{Notation}

%%%%%

\begin{abstract}
We prove an estimate on the difference of Maslov indices relative to the choice of two distinct
reference Lagrangians of a continuous path in the Lagrangian Grassmannian of a symplectic space.
We discuss some applications to the study of conjugate and focal points along a geodesic in
a semi-Riemannian manifold.
\end{abstract}

\maketitle
%\tableofcontents
%%%%%%%%%%%%%%%%%%%%%%%%%%%%%%%%%%%%%%%%%%%%%%%%%%%%%%%%%%%%%%%%
\begin{section}{Introduction}

Classical comparison theorems for conjugate and focal points in Riemannian or causal
Lorentzian geometry require curvature assumptions, or Morse theory (see \cite{Amb, EscOsu, Gal, Kup0, Kup1}). When passing to
the general semi-Riemannian world this approach does not work. Namely, the curvature is
never bounded (see \cite{DaNo}) and the index form has always infinite Morse index. In addition, it is
well known that singularities of the semi-Riemannian exponential map may accumulate
along a geodesic (see \cite{fechado}), and there is no hope to formulate a meaningful
comparison theorem using assumptions on the number of conjugate or focal points.

There are several good indications that a suitable substitute of the notion of size of the set of conjugate
or focal points along a semi-Riemannian geodesic is given by the Maslov index. This is a symplectic
integer valued invariant associated to the Jacobi equation, or more generally to the linearized
Hamilton equations along the solution of a Hamiltonian system. This number replaces the Morse
index of the index form, which in the general semi-Riemannian case is always infinite,
and in some \emph{nondegenerate} case it is a sort of algebraic count of the conjugate
points. In the Riemannian or causal Lorentzian case, the Maslov index of a geodesic relative to some fixed
Lagrangian coincides with the number of conjugate (or focal) points counted with multiplicity.
The exponential map is not locally injective around nondegenerate conjugate points (see \cite{War}),
or more generally around conjugate points whose contribution to the Maslov index is non zero
(see \cite{geobif}).

Inspired by a recent article by A. Lytchak \cite{Lyt},
in this paper we prove an estimate on the difference between Maslov indices
(Proposition~\ref{thm:stimacentrale}), and we apply this
estimate to obtain a number of results that are the semi-Riemannian analogue of the standard
comparison theorems in Riemannian geometry (Section~\ref{sec:comparison}).
These results relate the existence and the multiplicity
of conjugate and focal points with the values of Maslov indices naturally associated to a given geodesic.
It is very interesting to observe that Riemannian versions of the results proved in the present paper,
which are mostly well known, are obtained here with a proof that appears to be significantly
more elementary than the classical proof using Morse theory.

The paper is organized as follows. In Section~\ref{sec:preliminaries} we recall a few basic facts
on the geometry of the Lagrangian Grassmannian $\Lambda$ of a symplectic space $(V,\omega)$,
and on the notion of Maslov index for continuous paths in $\Lambda$. We use a generalized notion
of Maslov index, which applies to paths with arbitrary endpoints; note that, for paths with endpoints
on the Maslov cycle, there are several conventions regarding the contribution of the endpoints.
Here we adopt a convention slightly different from that in \cite{RobSal1},
(see \eqref{eq:primaopzione}, \eqref{eq:relMaslovneg} and \eqref{eq:secondaopzione}),
which is better suited for our purposes.

Section~\ref{sec:estimate} contains the estimate \eqref{eq:best}
on the difference of Maslov indices relatively
to the choice of two arbitrarily fixed reference Lagrangians $L_0$ and $L_1$.
Using the canonical atlas of charts of the Grassmannian Lagrangian and the transtition map
\eqref{eq:transitionmap}, the proof is
reduced to studying the index of perturbations of symmetric bilinear forms (Lemma~\ref{thm:stimacoindice},
Corollary~\ref{thm:stimadifferenza}). Several analogous estimates (\eqref{eq:best2}, \eqref{eq:best3})
are obtained using the properties \eqref{eq:relMaslovneg} and \eqref{eq:symmetriesHormander} of H\"ormander's index.

Applications to the study of conjugate and focal points along semi-Riemannian geodesics
are discussed in Section~\ref{sec:comparison}. In Subsection~\ref{sub:geolagr} we describe how
to obtain Lagrangian paths out of the flow of the Jacobi equation along a geodesic $\gamma:[a,b]\to M$
and an initial nondegenerate submanifold $\mathcal P$ of a semi-Riemannian manifold $(M,g)$.
In Lemma~\ref{thm:charLagrinitialsubmanifold} we give a characterization of which Lagrangian
subspaces of the symplectic space $T_{\gamma(a)}M\oplus T_{\gamma(a)}M$ arise from an initial
submanifold construction. The comparison results are proved in Subsection~\ref{sub:comparisonresult};
they include comparison between conjugate and focal points, as well as comparison between conjugate points
relative to distinct initial endpoints.
We conclude the paper in Section~\ref{sec:final} with a few final remarks concerning the question
of nondegeneracy of conjugate and focal points.

\end{section}
\begin{section}{Preliminaries}
\label{sec:preliminaries}
\subsection{The Lagrangian Grassmannian}
Let us consider a symplectic space $(V,\omega)$, with $\Dim(V)=2n$;
we will denote by $\Spl(V,\omega)$ the \emph{symplectic group} of $(V,\omega)$,
which is the closed Lie subgroup of $\GL(V)$ consisting of all isomorphisms
that preserve $\omega$.
A subspace $X\subset V$ is \emph{isotropic} if the restriction of $\omega$ to $X\times X$ vanishes
identically; an $n$-dimensional (i.e., maximal) isotropic subspace $L$ of $V$ is called
a \emph{Lagrangian subspace}. We denote by $\Lambda$ the Lagrangian Grassmannian of $(V,\omega)$,
which is the collection of all Lagrangian subspaces of $(V,\omega)$, and is a compact differentiable
manifold of dimension $\frac12n(n+1)$. A real-analytic atlas of charts on $\Lambda$ is given
as follows. Given a Lagrangian decomposition $(L_0,L_1)$ of $V$, i.e., $L_0,L_1\in\Lambda$ are transverse
Lagrangians, so that $V=L_0\oplus L_1$, then denote by $\Lambda^0(L_1)$ the open and dense subset
of $\Lambda$ consisting of all Lagrangians $L$ transverse to $L_1$. A diffeomorphism
$\varphi_{L_0,L_1}$ from $\Lambda^0(L_1)$ to the vector space $\Bsym(L_0)$ of all symmetric
bilinear forms on $L_0$ is defined by $\varphi_{L_0,L_1}(L)=\omega(T\cdot,\cdot)\vert_{L_0\times L_0}$,
where $T:L_0\to L_1$ is the unique linear map whose graph in $L_0\oplus L_1=V$ is $L$.
The kernel of $\varphi_{L_0,L_1}(L)$ is the space $L\cap L_0$.

We will need the following expression for the transition map $\varphi_{L_1,L}\circ\varphi_{L_0,L}^{-1}$,
where $L_0,L_1,L\in\Lambda$ are three Lagrangians such that $L\cap L_0=L\cap L_1=\{0\}$.
Note that the two charts $\varphi_{L_0,L}$ and $\varphi_{L_1,L}$ have the same domain.
If $\eta:L_1\to L_0$ denotes the isomorphism defined as the restriction to $L_1$ of the projection
$L\oplus L_0\to L_0$, then for all $B\in\Bsym(L_0)$ the following formula holds (see for instance  \cite[Lemma~2.5.4]{notasXIescola}):
\begin{equation}\label{eq:transitionmap}
\varphi_{L_1,L}\circ\varphi_{L_0,L}^{-1}(B)=\eta^*B+\varphi_{L_1,L}(L_0),
\end{equation}
where $\eta^*$ is the pull-back by $\eta$.

If $(L_0,L_1)$ is a Lagrangian decomposition of $V$, there exists a bijection between $\Lambda$ and
the set of pairs $(P,S)$, where $P\subset L_1$ is a subspace and $S:P\times P\to\R$ is a symmetric bilinear
form on $P$ (see \cite[Exercise 1.11]{notasXIescola}). More precisely, to each pair $(P,S)$ one associates the Lagrangian subspace $L_{P,S}$ defined
by:
\begin{equation}\label{eq:LagrPS}
L_{P,S}=\Big\{v+w:v\in P,\ w\in L_0,\ \omega(w,\cdot)\vert_P+S(v,\cdot)=0\Big\}.
\end{equation}
\subsection{Maslov index}
Let us recall a few notions related to symmetric bilinear forms.
Given a symmetric bilinear form $B$ on a (finite dimensional) real vector space $W$,
the \emph{index} of $B$ is defined to be the dimension of a maximal subspace of $W$ on which
$B$ is negative definite. The coindex of $B$ is the index of $-B$, and the {\em signature}
of $B$, denoted by $\segnatura(B)$ is defined to be the difference coindex minus index.

We will now recall briefly the notion of Maslov index for a continuous path $\ell:[a,b]\to\Lambda$.
For a fixed Lagrangian $L_0\in\Lambda$, the \emph{$L_0$-Maslov index }$\mu_{L_0}(\ell)$ of
$\ell$ is the integer characterized by the following properties:
\begin{itemize}
\item[(a)] $\mu_{L_0}$ is fixed-endpoint homotopy invariant;
\item[(b)] $\mu_{L_0}$ is additive by concatenation;
\item[(c)] if $\ell\big([a,b])\subset\Lambda^0(L_1)$ for some Lagrangian $L_1$ transverse
to $L_0$, then
\begin{equation}\label{eq:primaopzione}
\mu_{L_0}(\ell)=n_+\big[\varphi_{L_0,L_1}\big(\ell(b)\big)\big]-
n_+\big[\varphi_{L_0,L_1}\big(\ell(a)\big)\big],
\end{equation}
\end{itemize}
(see \cite{GPP04} for a similar discussion).
Let us denote by $\mu_{L_0}^-$ the $L_0$-Maslov index function relatively to the opposite
symplectic form $-\omega$ on $V$. The relation between the functions $\mu_{L_0}$ and
$\mu^-_{L_0}$ is given by the following identity:
\begin{equation}\label{eq:relMaslovneg}
\mu_{L_0}^-(\ell)=-\mu_{L_0}(\ell)+\Dim\big(\ell(a)\cap L_0\big)-\Dim\big(\ell(b)\cap L_0\big),
\end{equation}
for every continuous path $\ell:[a,b]\to\Lambda$.

Let us emphasize that, for curves $\ell$ whose endpoints are not transverse to $L_0$, there
are several conventions as to the contribution to the Maslov index of the endpoints.
For instance, the definition of $L_0$-Maslov index $\bar\mu_{L_0}$ in \cite{RobSal1} is\footnote{%
With such convention, the
Maslov index changes sign when one takes the opposite symplectic form.} obtained by replacing
\eqref{eq:primaopzione} with:
\begin{equation}\label{eq:secondaopzione}
\bar\mu_{L_0}(\ell)=\tfrac12\segnatura\big[\varphi_{L_0,L_1}\big(\ell(b)\big)\big]-
\tfrac12\segnatura\big[\varphi_{L_0,L_1}\big(\ell(a)\big)\big],
\end{equation}
in which case the Maslov index takes values in $\frac12\Z$.

\smallskip

Given any continuous path $\ell:[a,b]\to\Lambda$ and any two Lagrangians $L_0,L_0'\in\Lambda$,
the difference $\mu_{L_0}(\ell)-\mu_{L_0'}(\ell)$ depends only on $L_0$, $L_0'$ and
the endpoints $\ell(a)$ and $\ell(b)$ of $\ell$. This quantity will be denoted by
$\mathfrak q\big(L_0,L_0';\ell(a),\ell(b)\big)$, and it coincides  (up to some factor which is
irrelevant here) with the so called \emph{H\"ormander index} (see \cite{Hor}).
The H\"ormander index satisfies certain symmetries; we will need the following:
\begin{equation}\label{eq:symmetriesHormander}
\phantom{,\quad\forall\,L_0,L_1,L_0'}
\mathfrak q(L_0,L_1;L_0',L_1')=-\mathfrak q(L_0',L_1';L_0,L_1),\quad\forall\,L_0,L_1,L_0',L_1'\in\Lambda.
\end{equation}
The quantity:
\begin{equation}\label{eq:defKashiwara}
\tau(L_0,L_1,L_2)=\mathfrak q(L_0,L_1;L_2,L_0)=-\mathfrak q(L_0,L_1;L_0,L_2)
\end{equation}
coincides (again up to some factor) with the \emph{Kashiwara index} (see \cite{LioVer}). The Kashiwara index
function determines completely the H\"ormander index, by the identity:
\begin{equation}\label{eq:kashiwarahormander}
\mathfrak q(L_0,L_1;L_0',L_1')=\tau(L_0,L_1,L_0')-\tau(L_0,L_1,L_1'),\quad\forall\,L_0,L_1,L_0',L_1'\in\Lambda,
\end{equation}
which is easily proved using the concatenation additivity property of the Maslov index.

\end{section}

\begin{section}{An estimate on the difference of Maslov indices}
\label{sec:estimate}
Our analysis is based on the following elementary result:
\begin{lem}\label{thm:stimacoindice}
Let $B$ and $C$ be symmetric bilinear forms on a (finite dimensional) real vector space $V$.
Then:
\[-n_-(C)\le n_+(B+C)-n_+(B)\le n_+(C).\]
\end{lem}
\begin{proof}
It suffices to prove the inequality $n_+(B+C)-n_+(B)\le n_+(C)$; if this holds for every $B$ and $C$, replacing
$C$ with $-C$ and $B$ with $B+C$ will yield the other inequality $-n_-(C)\le n_+(B+C)-n_+(B)$.
Choose $W\subset V$ a maximal subspace of $V$ on which $B+C$ is positive definite,
so that $\Dim(W)=n_+(B+C)$, and write $W=W_+\oplus W_-$, where $B\vert_{W_+\times W_+}$ is
positive definite and $B\vert_{W_-\times W_-}$ is negative semi-definite.
Since $B+C$ is positive definite on $W$, it follows that $C\vert_{W_-\times W_-}$ must be positive definite,
so that $n_+\big(C\vert_{W\times W}\big)\ge\Dim(W_-)$.
Then:
\begin{multline*}
n_+(B+C)=\Dim(W)=\Dim(W_-)+\Dim(W_+)\\\le n_+\big(C\vert_{W\times W}\big)+n_+\big(B\vert_{W\times W}\big)
\le n_+(B)+n_+(C).\qedhere
\end{multline*}
\end{proof}
\begin{cor}\label{thm:stimadifferenza}
Given a fixed symmetric bilinear form $C$ on $V$, then for all $B_1,B_2\in\Bsym(V)$:
\[\big\vert n_+(B_1)-n_+(B_2)-n_+(B_1+C)+n_+(B_2+C)\big\vert\le n_-(C)+n_+(C).\qed\]
\end{cor}
\begin{prop}\label{thm:stimacentrale}
Given any continuous curve $\ell:[a,b]\to\Lambda$ and any pair $L_0,L_1\in\Lambda$ of Lagrangians,
then:
\begin{equation}\label{eq:best}
\big\vert\mu_{L_0}(\ell)-\mu_{L_1}(\ell)\big\vert\le n-\Dim(L_0\cap L_1).
\end{equation}
\end{prop}
\begin{proof}
Since the quantity $\mu_{L_0}(\ell)-\mu_{L_1}(\ell)$ depends only on the endpoints
$\ell(a)$ and $\ell(b)$, we can assume the existence of a Lagrangian $L\in\Lambda^0(L_0)\cap\Lambda^0(L_1)$
such that $\ell(t)\in\Lambda^0(L)$ for all $t\in[a,b]$. Namely, one can choose $L\in\Lambda^0(L_0)\cap
\Lambda^0(L_1)\cap\Lambda^0\big(\ell(a)\big)\cap\Lambda^0\big(\ell(b)\big)$ (these are dense opens subsets
of $\Lambda$, hence their intersection is non empty!), and replace $\ell$ by any continuous curve
in $\Lambda^0(L)$ from $\ell(a)$ to $\ell(b)$.

Once we are in this situation, then the Maslov indices of $\ell$ are given by:
\[\mu_{L_0}(\ell)=n_+\big[\varphi_{L_0,L}\big(\ell(b)\big)\big]-n_+\big[\varphi_{L_0,L}\big(\ell(a)\big)\big],\]
\[\mu_{L_1}(\ell)=n_+\big[\varphi_{L_1,L}\big(\ell(b)\big)\big]-n_+\big[\varphi_{L_1,L}\big(\ell(a)\big)\big].\]
Now consider the isomorphism $\eta:L_1\to L_0$ obtained as the restriction to $L_1$ of the projection
$L\oplus L_0\to L_0$; using formula \eqref{eq:transitionmap} of transition function for the charts $\varphi_{L_0,L}$ and $\varphi_{L_1,L}$, for all $\alpha\in\Lambda^0(L)$ we have:
\[\varphi_{L_1,L}(\alpha)=\eta^*\big(\varphi_{L_0,L}(\alpha)+\eta_*\varphi_{L_1,L}(L_0)\big),\]
and so:
\[n_+\big(\varphi_{L_1,L}(\alpha)\big)=n_+\big(\varphi_{L_0,L}(\alpha)+C\big),\]
where:
\[C=\eta_*\varphi_{L_1,L}(L_0)\]
does not depend on $\alpha$. Note that:
\[n_+(C)+n_-(C)=n-\Dim\big(\Ker(C)\big)=n-\Dim(L_0\cap L_1).\]
Inequality \eqref{eq:best} is obtained easily from Corollary~\ref{thm:stimadifferenza}
by setting $B_1=\varphi_{L_0,L}\big(\ell(b)\big)$ and $B_2=\varphi_{L_0,L}\big(\ell(a)\big)$.
\end{proof}
Using the symmetry property \eqref{eq:symmetriesHormander} of H\"ormander index, we also
get the following estimate:
\begin{cor}\label{eq:stimalaternativa}
Given any continuous curve $\ell:[a,b]\to\Lambda$ and any pair $L_0,L_1\in\Lambda$ of Lagrangians,
then:
\begin{equation}\label{eq:best2}
\big\vert\mu_{L_0}(\ell)-\mu_{L_1}(\ell)\big\vert\le n-\Dim\big(\ell(a)\cap\ell(b)\big).\qed
\end{equation}
\end{cor}
Moreover, changing the sign of the symplectic form and using \eqref{eq:relMaslovneg}, one obtains easily the following inequalities:
\begin{gather}
\notag\begin{aligned}
\big\vert\,\mu_{L_0}(\ell)-\mu_{L_1}(\ell)-\Dim\big(\ell(a)\cap L_0\big)&+\Dim\big(\ell(a)\cap L_1\big)\\+\Dim\big(\ell(b)\cap L_0\big)&
-\Dim\big(\ell(b)\cap L_1\big)\big\vert\le n-\Dim(L_0\cap L_1),
\end{aligned}
\\ \label{eq:best3}
\begin{aligned}
\big\vert\,\mu_{L_0}(\ell)-\mu_{L_1}(\ell)-\Dim\big(\ell(a)\cap L_0\big)&+\Dim\big(\ell(a)\cap L_1\big)\\+\Dim\big(\ell(b)\cap L_0\big)&
-\Dim\big(\ell(b)\cap L_1\big)\big\vert\le n-\Dim\big(\ell(a)\cap \ell(b)\big).
\end{aligned}
\end{gather}
\end{section}

\begin{section}{Comparison results for conjugate and focal points}
\label{sec:comparison}
\subsection{Geodesics and Lagrangian paths}
\label{sub:geolagr}
Let us now look more specifically at curves of Lagrangians arising from the Jacobi equation
along a semi-Riemannian geodesic. Let $(M,g)$ be a semi-Riemannian manifold of dimension $n$, $\nabla$ the
covariant derivative of the Levi--Civita connection of $g$, with curvature tensor
chosen with the sign convention $R(X,Y)=[\nabla_X,\nabla_Y]-\nabla_{[X,Y]}$.
We will assume throughout the section that $\gamma:[a,b]\to M$ is a given geodesic
in $M$; when needed, we will also consider extensions of $\gamma$ to a larger interval $[a',b']\supset[a,b]$.
The Jacobi equation along $\gamma$ is given by $\Ddt^2V-R(\dot\gamma,V)\dot\gamma=0$.
Consider the flow of the Jacobi equation,
which is the family of isomorphisms \[\Phi_t:T_{\gamma(a)}M\oplus T_{\gamma(a)}M\longrightarrow T_{\gamma(t)}M\oplus T_{\gamma(t)}M,\] $t\in[a,b]$, defined by $\Phi_t(v,w)=\big(J_{v,w}(t),\Ddt J_{v,w}(t)\big)$, where
$J_{v,w}$ is the unique Jacobi field along $\gamma$ satisfying $J(a)=v$ and $\Ddt J(a)=w$.
Consider the symplectic form $\omega$ on the space $V=T_{\gamma(a)}M\oplus T_{\gamma(a)}M$
given by $\omega\big(v_1,w_1),(v_2,w_2)\big)=g(v_2,w_1)-g(v_1,w_2)$. For all $t\in[a,b]$,
define $L_0^t=\{0\}\oplus T_{\gamma(t)}M\subset T_{\gamma(t)}M\oplus T_{\gamma(t)}M$ and set
$\ell(t)=\Phi_t^{-1}(L_0^t)$. An immediate calculation shows that $\ell(t)$ is a Lagrangian
subspace of $(V,\omega)$, and we obtain in this way a smooth curve $\ell:[a,b]\to\Lambda(V,\omega)$.
Note that:
\begin{equation}\label{eq:valoreiniziale}
\ell(a)=L^a_0=: L_0.
\end{equation}
Now, consider a smooth connected submanifold $\mathcal P\subset M$, with $\gamma(a)\in\mathcal P$ and\footnote{%
In this section, the symbol $\perp$ will denote orthogonality with respect to the semi-Riemannian
metric $g$.}
$\dot\gamma(a)\in T_{\gamma(a)}\mathcal P^\perp$; let us also assume that $\mathcal P$
is nondegenerate at $\gamma(a)$, meaning that the restriction of the metric $g$ to
$T_{\gamma(a)}\mathcal P$ is nondegenerate.
We will denote by $n_-(g,\mathcal P)$ and $n_+(g,\mathcal P)$ respectively the index and the coindex
of the restriction of $g$ to $\mathcal P$, so that $n_-(g,\mathcal P)+n_+(g,\mathcal P)=\Dim(\mathcal P)$.
Let $S$
be the second fundamental form of $\mathcal P$ at $\gamma(a)$ in the normal direction $\dot\gamma(a)$,
seen as a $g$-symmetric operator $S:T_{\gamma(a)}\mathcal P\to T_{\gamma(a)}\mathcal P$,
and consider the subspace $L_{\mathcal P}\subset V$ defined by:
\[L_{\mathcal P}=\Big\{(v,w)\in T_{\gamma(a)}M\oplus T_{\gamma(a)}M:v\in T_{\gamma(a)}\mathcal P,\ w+S(v)\in T_{\gamma(a)}\mathcal P^\perp\Big\},\]
which is precisely the construction of Lagrangian subspaces described abstractly in \eqref{eq:LagrPS}.
If $\pi_1:T_{\gamma(a)}M\oplus T_{\gamma(a)}M\to T_{\gamma(a)}M$ is the projection onto the first
summand, then $\pi_1(L_{\mathcal P})=T_{\gamma(a)}\mathcal P$ is orthogonal to $\dot\gamma(a)$.
Conversely:
\begin{lem}\label{thm:charLagrinitialsubmanifold}
Let $L\subset T_{\gamma(a)}M\oplus T_{\gamma(a)}M$ be a Lagrangian subspace, and
assume that $P=\pi_1(L)$ is orthogonal to $\dot\gamma(a)$. Then, there exists a smooth submanifold
$\mathcal P$ orthogonal to $\dot\gamma(a)$ such that $L=L_{\mathcal P}$.
\end{lem}
\begin{proof}
Consider the Lagrangian decomposition $(L_0,L_1)$ of $T_{\gamma(a)}M\oplus T_{\gamma(a)}M$ given
by $L_0=\{0\}\oplus T_{\gamma(a)}M$ and $L_1=T_{\gamma(a)}M\oplus\{0\}$; then there exists
a symmetric bilinear form $S:P\times P\to R$ such that $L=L_{P,S}$ as in \eqref{eq:LagrPS}.
Let $\mathcal P_0\subset T_{\gamma(a)}M$ be the submanifold
given by the graph of the function $P\ni x\mapsto \frac12S(x,x)\dot\gamma(a)\in P^\perp$.
The desired submanifold $\mathcal P$ is obtained by taking the exponential of a small open neighborhood
of $0$ in $\mathcal P_0$. It is easily seen that the tangent space to $\mathcal P_0$ at $0$
is $P$, and since $\mathrm d\exp_{\gamma(a)}(0)$ is the identity, $T_{\gamma(a)}\mathcal P=P$.
Moreover, using the fact that the Christoffel symbols of the chart $\exp_{\gamma(a)}$ vanish at $0$,
it is easily seen that the second fundamental form of $\mathcal P$ at $\gamma(a)$ in the normal
direction $\dot\gamma(a)$ is $S$.
\end{proof}

Let us also consider the space $L_0=\{0\}\oplus T_{\gamma(a)}M$, which corresponds to the Lagrangian
associated to the trivial initial submanifold $\mathcal P=\{\gamma(a)\}$.
Then, an instant $t\in\left]a,b\right]$ is $\mathcal P$-focal along $\gamma$ if and only if
$\ell(t)\cap L_{\mathcal P}\ne\{0\}$, and the dimension of this intersection equals the multiplicity
of $t$ as a $\mathcal P$-focal instant. In particular, $t$ is a conjugate instant, i.e., $\gamma(t)$
 is conjugate to $\gamma(a)$ along $\gamma$, if $\ell(t)\cap L_0\ne\{0\}$.
 Note that:
 \begin{equation}\label{eq:L0capLP}
 L_0\cap L_{\mathcal P}=\{0\}\oplus T_{\gamma(a)}\mathcal P^\perp,
 \end{equation}
 thus:
 \begin{equation}\label{eq:dimL0capLP}
 \Dim(L_0\cap L_{\mathcal P})=\codim(\mathcal P).
 \end{equation}

 For all $t\in\left]a,b\right]$, consider the space
 \[A_{\mathcal P}[t]=\Big\{\Ddt J(t): \text{$J$ is a $\mathcal P$-Jacobi field along $\gamma$ with}\ J(t)=0\Big\},\]
 while for $t=a$ we set:
 \[A_{\mathcal P}[a]=T_{\gamma(a)}\mathcal P^\perp;\]
 note that $\Dim\big(A_{\mathcal P}[t]\big)=\Dim\big(\ell(t)\cap L_{\mathcal P}\big)$.
 When the initial submanifold is just a point, we will use the following notation:
 \begin{equation}\label{eq:defA0t}
 \begin{gathered}
 A_0[t]=\Big\{\Ddt J(t): \text{$J$ is a Jacobi field along $\gamma$ with $J(a)=0$ and $J(t)=0$}\Big\},\\\
 A_0[a]=T_{\gamma(a)}M.
 \end{gathered}
 \end{equation}
 It is well known that focal or conjugate points along a semi-Riemannian geodesic may accumulate (see \cite{fechado}),
 however, \emph{nondegenerate} conjugate or focal points are isolated. A $\mathcal P$-focal point
 $\gamma(t)$ along $\gamma$ is nondegenerate when the restriction of the metric $g$ to the space
 $A_{\mathcal P}[t]$ is nondegenerate.
 This is always the case when $g$ is positive definite (i.e., Riemannian), or
 if $g$ has index $1$ (i.e., Lorentzian) and $\gamma$ is either timelike or lightlike.
 Also, the initial endpoint $\gamma(a)$ which is
 always $\mathcal P$-focal of multiplicity equal to the codimension of $\mathcal P$, is always isolated.

 For all $t\in[a,b]$, let us denote by $n_-(g,\mathcal P,t)$, $n_+(g,\mathcal P,t)$ and
 $\sigma(g,\mathcal P,t)$ respectively the index, the coindex and the signature of the restriction
 of $g$ to $A_{\mathcal P}[t]$.
 Given a nondegenerate $\mathcal P$-focal point $\gamma(t)$ along $\gamma$, with $t\in\left]a,b\right[$, then
 $t$ is an isolated instant of nontransversality of the Lagrangians $\ell(t)$ and $L_{\mathcal P}$.
  Its contribution to the Maslov index $\mu_{L_{\mathcal P}}(\ell)$, i.e.,
  $\mu_{L_{\mathcal P}}(\ell\vert_{[t-\varepsilon,t+\varepsilon]})$ with $\varepsilon>0$ sufficiently small, is given by
 the integer $\sigma(g,\mathcal P,t)$.
 The contribution of the initial point  to the Maslov index $\mu_{L_{\mathcal P}}(\ell)$, which as observed is always nondegenerate,
 is given by $n_+(g,\mathcal P,a)$:
 \begin{equation}\label{eq:maslovinizio}
 \mu_{L_\mathcal P}\big(\ell\vert_{[a,a+\varepsilon]}\big)=n_+(g,\mathcal P,a)=n_+(g)-n_+(g,\mathcal P).
 \end{equation}
 In particular:
 \begin{equation}\label{eq:maslovinizio2}
 \mu_{L_0}\big(\ell\vert_{[a,a+\varepsilon]}\big)=n_+(g).
 \end{equation}

 Moreover, if $\gamma(b)$ is a nondegenerate $\mathcal P$-focal
 point along $\gamma$, then its contribution to the Maslov index  $\mu_{L_{\mathcal P}}(\ell)$ is equal
 to $-n_-(g,\mathcal P,b)$.
 Thus, when $g$ is Riemannian the Maslov index $\mu_{L_{\mathcal P}}\big(\ell\vert_{[a+\varepsilon,b]}\big)$
 is the number of $\mathcal P$-focal points along $\gamma\vert_{\left[a,b\right[}$
 counted with multiplicity. The same holds when $g$ is Lorentzian (i.e., index equal to $1$)
 and $\gamma$ is timelike. More generally, if all $\mathcal P$-focal points along $\gamma$
 are nondegenerate, the Maslov index $\mu_{L_{\mathcal P}}(\ell)$ is given by the finite sum:
 \[\mu_{\mathcal P}(\ell)=n_+(g)-n_+(g,\mathcal P)+\sum_{t\in\left]a,b\right[}\sigma(g,\mathcal P,t)-n_-(g,\mathcal P,b).\]
 All this follows easily from the following elementary result:
 \begin{lem}\label{thm:elementary}
 Let $B:I\to\Bsym(V)$ be a $C^1$-curve of symmetric bilinear forms on a real vector space $V$.
 Assume that $t_0\in I$ is a degeneracy instant, and denote by $B_0$ the restriction to $\Ker\big(B(t_0)\big)$
 of the derivative $B'(t_0)$. If $B_0$ is nondegenerate, then $t_0$ is an isolated
 degeneracy instant, and for $\varepsilon>0$ sufficiently small:
 \[n_+\big(B(t_0+\varepsilon)\big)-n_+\big(B(t_0)\big)=n_+(B_0),\quad n_+\big(B(t_0)\big)-n_+\big(B(t_0-\varepsilon)\big)=-n_-(B_0).\qed\]
 \end{lem}
 Lemma~\ref{thm:elementary} is employed in order to compute the Maslov index $\mu_{L_{\mathcal P}}$ as follows.
 Given a $\mathcal P$-focal instant $t_0\in[a,b]$
 and a Lagrangian $L_1$ transversal to both $L_{\mathcal P}$ and $\ell(t_0)$,
 then consider the smooth path $t\mapsto\varphi_{L_{\mathcal P},L_1}\big(\ell(t)\big)$
 of symmetric bilinear forms on $L_{\mathcal P}$.  The kernel of $B(t_0)$ is identified with the
 space $A_{\mathcal P}[t_0]$, and the restriction of the derivative $B'(t_0)$ to $\Ker\big(B(t_0)\big)$
 with the restriction of the metric $g$ to $A_{\mathcal P}[t_0]$ (see for instance \cite{pacific}).
 \subsection{Comparison results}
 \label{sub:comparisonresult}
 Having this in mind, let us now prove some comparison results for conjugate and focal instants.
  \begin{prop}\label{thm:stimamaslovgeo}
 Given any interval $[\alpha,\beta]\subset[a,b]$:
 \begin{equation}\label{eq:stimemaslovgeo}
 \left\vert\,\mu_{L_0}\big(\ell\vert_{[\alpha,\beta]}\big)-\mu_{L_{\mathcal P}}\big(\ell\vert_{[\alpha,\beta]}\big)\right\vert
 \le\Dim(\mathcal P).
 \end{equation}
 \end{prop}
 \begin{proof}
 It follows readily from Proposition~\ref{thm:stimacentrale} and \eqref{eq:dimL0capLP}.
 \end{proof}
 In particular, we have the following result concerning the existence of conjugate or focal
 instant along an arbitrary portion of a geodesic:
 \begin{cor}
Given any interval $[\alpha,\beta]\subset\left]a,b\right]$:
\begin{itemize}
\item if $\big\vert\,\mu_{L_0}\big(\ell\vert_{[\alpha,\beta]}\big)\big\vert>\Dim(\mathcal P)$, then there is
at least one $\mathcal P$-focal instant in $[\alpha,\beta]$;
\item if $\big\vert\,\mu_{L_{\mathcal P}}\big(\ell\vert_{[\alpha,\beta]}\big)\big\vert>\Dim(\mathcal P)$, then there is
at least one conjugate instant in $[\alpha,\beta]$.
\end{itemize}
 \end{cor}
 \begin{proof}
 By Proposition~\ref{thm:stimamaslovgeo}, if $\big\vert\,\mu_{L_0}\big(\ell\vert_{[\alpha,\beta]}\big)\big\vert>\Dim(\mathcal P)$
 then $\big\vert\,\mu_{L_{\mathcal P}}\big(\ell\vert_{[\alpha,\beta]}\big)\big\vert>0$.
 Since $a\not\in[\alpha,\beta]$, this implies that there is a $\mathcal P$-focal instant in $[\alpha,\beta]$.
 The second statement is totally analogous.
 \end{proof}
 On the other hand, the absence of conjugate (focal) instants gives an upper bound on the number
 of focal (conjugate) instants:
 \begin{prop}
 If $\gamma$ has no conjugate instant, then for every interval $[\alpha,\beta]\subset\left]a,b\right]$,
 $\left\vert\,\mu_{L_{\mathcal P}}\big(\ell\vert_{[\alpha,\beta]}\big)\right\vert\le\Dim(\mathcal P)$.
 Similarly, if $\gamma$ has no $\mathcal P$-focal instant, then
 $\left\vert\,\mu_{L_{0}}\big(\ell\vert_{[\alpha,\beta]}\big)\right\vert\le\Dim(\mathcal P)$.
 \end{prop}
 \begin{proof}
 If $\gamma$ has no conjugate (resp., $\mathcal P$-focal)
 instant, then the Maslov index $\mu_{L_{0}}\big(\ell\vert_{[\alpha,\beta]}\big)=0$
 (resp., $\mu_{L_{\mathcal P}}\big(\ell\vert_{[\alpha,\beta]}\big)=0$) for
 all $[\alpha,\beta]\subset\left]a,b\right]$.
 \end{proof}
 All the above statements have a much more appealing version in the Riemannian or timelike Lorentzian
case, where the ``Maslov index'' can be replaced by the number of conjugate or focal instants.
In this situation, focal and conjugate instants are always nondegenerate and isolated,
and without using Morse theory one can prove nice comparison results of the following type:
\begin{cor}\label{thm:primafocali}
Assume that either $g$ is Riemannian or that $g$ is Lorentzian and $\gamma$ is timelike (in which
case $\mathcal P$ is necessarily a spacelike submanifold of $M$). Denote by $t_0$ and $t_{\mathcal P}$
the following instants:
\[t_0=\sup\Big\{t\in\left]a,b\right]:\ \text{there are no conjugate instants in $\left]a,t\right]$}\Big\},\]
\[t_{\mathcal P}=\sup\Big\{t\in\left]a,b\right]:\ \text{there are no $\mathcal P$-focal
instants in $\left]a,t\right]$}\Big\}.\]
Then, $t_{\mathcal P}\le t_0$, and if $t_{\mathcal P}=t_0$ then the multiplicity of $t_{\mathcal P}$
as a $\mathcal P$-focal point is greater than or equal to its multiplicity as a conjugate point.
\end{cor}
\begin{proof}
Assume $t_0<t_{\mathcal P}\le b$ and choose $t'\in\left]t_0,t_{\mathcal P}\right[$.
Since there are no $\mathcal P$-focal instants in $\left]a,t'\right]$ and $\mathcal P$ is spacelike, from \eqref{eq:maslovinizio} it follows that
$\mu_{L_{\mathcal P}}\big(\ell\vert_{[a,t']}\big)=\codim(\mathcal P)-n_-(g)$.
On the other hand, since $t_0$ is conjugate, $\mu_{L_0}\big(\ell\vert_{[a,t']}\big)\ge n_+(g)+1$,
hence:
\[\mu_{L_0}\big(\ell\vert_{[a,t']}\big)-\mu_{L_{\mathcal P}}\big(\ell\vert_{[a,t']}\big)\ge
n_+(g)+n_-(g)-\codim(\mathcal P)+1=\Dim(\mathcal P)+1,\]
contradicting \eqref{eq:stimemaslovgeo}.

Assume that $t_{\mathcal P}=t_0$ and that $t_{\mathcal P}$ is a $\mathcal P$-focal point.
By possibly extending the geodesic $\gamma$ to a slightly larger interval $[a,b']$ with
$b'>b$, we can assume the existence of $t'>t_{\mathcal P}$ with the property that
there are no conjugate or $\mathcal P$-focal instants in $\left]t_{\mathcal P},t'\right]$.
Then:
\[\mu_{L_0}\big(\ell\vert_{[a,t']}\big)=n_+(g)+\mathrm{mul}(t_{\mathcal P}),\]
where $\mathrm{mul}(t_{\mathcal P})$ is the (possibly null) multiplicity of $t_{\mathcal P}$ as
a conjugate instant.
Similarly:
\[\mu_{L_{\mathcal P}}\big(\ell\vert_{[a,t']}\big)=\codim(\mathcal P)-n_-(g)+\mathrm{mul}_{\mathcal P}(t_{\mathcal P}),\]
where $\mathrm{mul}_{\mathcal P}(t_{\mathcal P})$ is the multiplicity of $t_{\mathcal P}$ as
a $\mathcal P$-focal instant. Then:
\[\mu_{L_0}\big(\ell\vert_{[a,t']}\big)-\mu_{L_{\mathcal P}}\big(\ell\vert_{[a,t']}\big)=\Dim(\mathcal P)+
\mathrm{mul}(t_{\mathcal P})-\mathrm{mul}_{\mathcal P}(t_{\mathcal P})\]
which has to be less than or equal to $\Dim(\mathcal P)$, giving $\mathrm{mul}(t_{\mathcal P})
\ge\mathrm{mul}_{\mathcal P}(t_{\mathcal P})$.
\end{proof}

It is known that the result of Corollary~\ref{thm:primafocali} does not hold without the assumption
that the metric $g$ is positive definite or that $g$ is Lorentzian and $\gamma$ timelike.
A counterexample is exhibited by Kupeli in \cite{Kup1}, where the author constructs a
spacelike geodesic $\gamma$ orthogonal to a timelike submanifold $\mathcal P$ of a Lorentzian manifold,
with the property that $\gamma$ has conjugate points but no focal point.

In the following statements, $\varepsilon$ will denote a small positive number with the
property that there are no conjugate or $\mathcal P$-focal instants in $\left]a,a+\varepsilon\right]$.
\begin{prop}\label{thm:numefocmenonumconj}
The following inequalities hold:
\begin{equation*}
-n_-(g,{\mathcal{P}})\leq \mu_{L_{\mathcal{P}}}\big(\ell\vert_{[a+\varepsilon,b]}\big)
-\mu_{L_0}\big(\ell\vert_{[a+\varepsilon,b]}\big)\leq \dim \mathcal{P}.
\end{equation*}
\end{prop}
\begin{proof}
A straightforward consequence of formulas~\eqref{eq:maslovinizio}, \eqref{eq:maslovinizio2} and
\eqref{eq:stimemaslovgeo} applied on the intervals $[a,b]$ and $[a+\varepsilon,b]$.
\end{proof}
In particular, when $g$ is Riemannian, or $g$ is Lorentzian and $\gamma$ timelike,
Proposition~\ref{thm:numefocmenonumconj} says that the number of $\mathcal P$-focal points along $\gamma$ is greater than
or equal to the number of conjugate points along $\gamma$, and that their difference is less than or equal to  the
dimension of $\mathcal P$.
 \begin{cor}
 If $\mu_{L_0}\big(\ell\vert_{[a+\varepsilon,b]}\big)>n_-(g,\mathcal P)$ or
 $\mu_{L_0}\big(\ell\vert_{[a+\varepsilon,b]}\big)<-\Dim(\mathcal P)$, then there exists at least one
 $\mathcal P$-focal instant in $[a+\varepsilon,b]$.\qed
  \end{cor}

  \begin{cor}
 If $\mu_{L_{\mathcal P}}\big(\ell\vert_{[a+\varepsilon,b]}\big)<-n_-(g,\mathcal P)$ or
 $\mu_{L_{\mathcal P}}\big(\ell\vert_{[a+\varepsilon,b]}\big)>\Dim(\mathcal P)$, then there exists at least one
 conjugate instant in $[a+\varepsilon,b]$.\qed
  \end{cor}

For the following result we need to recall the definition of the space $A_0[t]$ given in \eqref{eq:defA0t};
we will denote by $n_+(g,t)$ and $n_-(g,t)$ respectively the coindex and the index of the restriction
of $g$ to $A_0[t]\times A_0[t]$ and $\mathrm{mul}(t_0)=\Dim\big(A_0[t_0]\big)$.

The estimate in Corollary~\ref{eq:stimalaternativa} can be used to obtain results of the following type:
\begin{cor}
If $t_0\in\left]a,b\right]$ is a conjugate
instant such that either:
$\mathrm{mul}(t_0)>n_-(g)-\mu_{L_0}\big(\ell\vert_{[a+\varepsilon,t_0]}\big)$ or $\mu_{L_0}\big(\ell\vert_{[a+\varepsilon,t_0]}\big)<-n_+(g)$, then
for every $a'<a$ there is an instant $t'\in[a,t_0]$ such that $\gamma(t')$ is conjugate to $\gamma(a)$ along $\gamma$.
\end{cor}
\begin{proof}
Consider the Lagrangian $L'\subset V$ given by:
\[L'=\Big\{(v,w)\in V:J_{v,w}(a')=0\Big\}.\]
If there were no instant
$t$ in $[a,t_0]$ with $\gamma(t)$ conjugate to $\gamma(a')$ along $\gamma$,
then $\mu_{L'}\big(\ell\vert_{[a,t_0]}\big)=\Dim\big(L'\cap\ell(a)\big)=\Dim\big(L'\cap\ell(t_0)\big)=0$.
By Corollary~\ref{eq:stimalaternativa} it would then be
\[\mu_{L_0}\big(\ell\vert_{[a,t_0]}\big)=\mu_{L_0}\big(\ell\vert_{[a,t_0]}\big)-\mu_{L'}\big(\ell\vert_{[a,t_0]}\big)\le n-\mathrm{mul}(t_0).\]
Using Eq. \eqref{eq:maslovinizio2}, we obtain a contradiction with the hypothesis of the corollary.
Moreover, using \eqref{eq:best3}, we have:
\begin{multline*}
\mathrm{mul}(t_0)-n=\Dim\big(\ell(a)\cap\ell(t_0)\big)-n\le\mu_{L_0}\big(\ell\vert_{[a,t_0]}\big)
-\Dim\big(\ell(a)\cap L_0\big)+\Dim\big(\ell(b)\cap L_0\big)\\=\mu_{L_0}\big(\ell\vert_{[a,t_0]}\big)-n
+\mathrm{mul}(t_0),
\end{multline*}
i.e.:
\[\mu_{L_0}\big(\ell\vert_{[a,t_0]}\big)\ge0,\]
which together with Eq. \eqref{eq:maslovinizio2} concludes the proof.
\end{proof}
When the  first conjugate point is nondegenerate, we can state a more precise result.
\begin{cor}\label{thm:cora'}
Let $t_0\in\left]a,b\right]$ be the first conjugate instant along $\gamma$, and assume that it is nondegenerate and $\mathrm{mul}(t_0)>n_-(g)+n_-(g,t_0)$. Then for every $a'<a$ 
there exists and instant $t'\in[a,t_0]$ such that
$\gamma(t')$ is conjugate to $\gamma(a')$ along $\gamma$. 
\end{cor}
Note that if $g$ is Riemannian, then $n_-(g)=n_-(g,t_0)=0$ and the result of Corollary~\ref{thm:cora'}
holds without any assumption of the multiplicity of $t_0$.

\end{section}

\begin{section}{Final remarks and conjectures}
\label{sec:final}
If the semi-Riemannian manifold $(M,g)$ is real-analytic,
then conjugate and focal points do not accumulate along a geodesic, and higher order formulas
for the contribution to the Maslov index of each conjugate and focal points are available (see \cite{lama}).
In this case, the statement of all the above results can be given in terms of the \emph{partial signatures}
of the conjugate and the focal points, which are a sort of generalized multiplicities.

It may also be worth observing that the nondegeneracy assumption for the conjugate and focal points
is stable by $C^3$-small perturbations of the metric, and \emph{generic}, although a precise
genericity statement seems a little involved to prove. We conjecture that, given a differentiable manifold
$M$ and a countable set $Z\subset TM$, then the set of semi-Riemannian  metrics $g$ on $M$ having a fixed
index and for which all the geodesics $\gamma:[0,1]\to M$ with $\dot\gamma(0)\in Z$ have only
conjugate points nondegenerate and of multiplicity equal to $1$ is generic.
In this situation, the comparison results proved in this paper would have a more explicit statement
in terms of number of conjugate and focal points.

A natural conjecture is also that in the case of \emph{stationary} Lorentzian metrics, all geodesics
have nondegenerate conjugate points whose contribution to the Maslov index is positive and equal to their
multiplicity. This fact has been proved in the case of left-invariant Lorentzian metrics on Lie groups
having dimension less than $6$ (see \cite{JaPi06a}) and, recently, using semi-Riemannian submersions (see \cite{CapJavPic}), also for spacelike
geodesics orthogonal to some timelike Killing vector field.
If this conjecture were true in full generality, one would have Riemannian-like comparison results also for spacelike
geodesics in stationary Lorentz manifolds.
\end{section}

\end{document}